# A remark on Sarnak's conjecture

Régis de la Bretèche & Gérald Tenenbaum

**Abstract.** We investigate Sarnak's conjecture on the Möbius function in the special case when the test function is the indicator of the set of integers for which a real additive function assumes a given value.

## 1. Introduction and statements of results

According to a general pseudo-randomness principle related to a famous conjecture of Chowla [1] and recently considered by Sarnak [7], the Möbius function $\mu$ does not correlate with any function $\xi$ of low complexity. In other words,

$$(1\cdot1) \qquad \sum_{n \leqslant x} \mu(n) \xi(n) = o\bigg( \sum_{n \leqslant x} |\xi(n)| \bigg) \qquad (x \to \infty).$$

There are many ways of constructing functions of low complexity. Sarnak and others use return times of sampling sequences of a dynamical system, which leads to a natural measure of the complexity. Here we propose to follow another path by selecting the test-function as the indicator of the set of those integers where a real additive function assumes a given value. It is known since Halász [5] that

$$(1\cdot2) \qquad Q(x;f) := \sup_{m \in \mathbb{R}} \sum_{\substack{n \leqslant x \\ f(n)=m}} 1 \ll \frac{x}{\sqrt{1+E(x)}}$$

where we have put

$$E(x) := \sum_{\substack{p \leqslant x \\ f(p) \neq 0}} \frac{1}{p}.$$

Here and in the sequel, the letter $p$ denotes a prime number.

The estimate (1·2) is known to be optimal in this generality since the two sides achieve the same order of magnitude when $f(n)$ is equal to the total number of prime factors of $n$, counted with or without multiplicity.

As a first investigation of the above described problem, we would like to show that

$$Q(x;f,\mu) := \sup_{m \in \mathbb{R}} \bigg| \sum_{\substack{n \leqslant x \\ f(n)=m}} \mu(n) \bigg|$$

is generically smaller than the right-hand side of (1·2). Of course we have to avoid the case when $f(p)$ is constant, for then $\mu(n)$ does not oscillate on the set of squarefree integers $n$ with $f(n) = m$. Therefore we seek an estimate which coincides with (1·2) when $f(p)$ is close to a constant and which has smaller order of magnitude otherwise.

When $f(p)$ is restricted to assume the values 0 or 1 only, we thus expect a significant improvement over (1·2) when

$$(1\cdot3) \qquad F(x) := \sum_{p \leqslant x} \frac{1-f(p)}{p}$$

is large. Indeed, in this simple case we obtain the following estimate.

**Theorem 1.1.** *Let $f$ denote a real additive arithmetic function such that $f(p) \in \{0,1\}$ for all $p$. Then, with the above notation and $c = (2\pi - 4)/(3\pi - 2) \approx 0.30751$, we have*

$$(1\cdot 4) \qquad Q(x; f, \mu) \ll \frac{x\{1 + F(x)\} e^{-cF(x)}}{\sqrt{1 + E(x)}}.$$

For simplicity, let us retain in the sequel the hypothesis $f(p) \in \{0,1\}$.[1] Under the assumption that $F(x)$, as defined in (1·3) above, grows sufficiently slowly, we may prove an estimate that is valid for each $m$ in a large range around the mean, and so may be stated in the exact frame of Sarnak's conjecture.

Let us denote by $N_m(x; f)$ the number of squarefree integers not exceeding $x$ such that $f(n) = m$. It follows from results of Halász [3], [4], and Sárközy [6] that, given any $\kappa \in\,]0, 1[$, we have

$$(1\cdot 5) \qquad N_m(x; f) \asymp x \frac{E(x)^m}{m!} e^{-E(x)} \qquad \bigl(\kappa E(x) \leqslant m \leqslant E(x)/\kappa\bigr).$$

Moreover, Halász announced (see [2], p. 312) the possibility to obtain, in the same range for $m$, an asymptotic formula for $N_m(x; f)$, a result which actually follows, as shown in [10], from a general effective mean value estimate for multiplicative functions established in the same work—see below.

This supports the hope to obtain an asymptotic formula for

$$N_m(x; f, \mu) := \sum_{\substack{n \leqslant x \\ f(n) = m}} \mu(n)$$

which directly compares to (1·5). In view of (1·1), we may assume with no loss of generality that $f$ is strongly additive. We obtain the following result. Here and in the sequel we let $\log_k$ denote the $k$-fold iterated logarithm.

**Theorem 1.2.** *Let $\kappa \in\,]0,1[$ and let $f$ denote a strongly additive function such that $f(p) \in \{0, 1\}$ for all primes $p$. Assume furthermore that*

$$(1\cdot 6) \qquad F(x) := \sum_{p \leqslant x} \frac{1 - f(p)}{p} \ll \log_3 x \qquad (x \to \infty)$$

$$(1\cdot 7) \qquad \sum_{\exp\{(\log x)/(\log_2 x)^D\} < p \leqslant y} \frac{\{1 - f(p)\} \log p}{p} \ll \frac{(\log y)}{(\log_2 x)^{c_0}} \qquad \bigl(x^{1/(\log_2 x)^D} < y \leqslant x\bigr)$$

*where $D$ and $c_0$ are positive constants. Provided $D$ is sufficiently large and uniformly in the range $\kappa E(x) \leqslant m \leqslant E(x)/\kappa$, we have*

$$(1\cdot 8) \qquad N_m(x; f, \mu) = (-1)^m N_m(x; f) \Bigl\{ \lambda_f e^{-2F(x)} + O\Bigl(\frac{1}{(\log_2 x)^b}\Bigr) \Bigr\},$$

*with*

$$(1\cdot 9) \qquad \lambda_f := \prod_{f(p)=0} \frac{1 - 1/p}{1 + 1/p} e^{2/p}, \quad b := \tfrac{1}{2} \min\{1, c_0 \kappa/(4 - \kappa)\}.$$

---

1. All our results could be straightforwardly adapted to case when $f(p)$ is restricted to a fixed, finite set, or even to a set of moderate size depending on $x$.

To fix ideas, note that a strongly additive function $f$ such that $f(p) \in \{0,1\}$ satisfies hypotheses (1·6) and (1·7) as soon as

$$\sum_{p \leqslant y} \{1 - f(p)\} \log p \ll \frac{y}{(\log_2 y)^{\max(1,c_0)}}.$$

The proof of Theorem 1.2 rests on the following recent result of the second author [10] (theorem 1.4), for the statement of which we introduce further notation. We let $\mathcal{M}(A, B)$ designate the class of those complex-valued multiplicative functions $g$ such that

(1·10) $$\max_p |g(p)| \leqslant A, \quad \sum_{p,\, \nu \geqslant 2} \frac{|g(p^\nu)| \log p^\nu}{p^\nu} \leqslant B,$$

and, for $\mathfrak{b} \in \mathbb{R}$, we write

(1·11) $$\beta_0 = \beta_0(\mathfrak{b}, A) := 1 - \frac{\sin(2\pi \mathfrak{b}/A)}{2\pi \mathfrak{b}/A}.$$

Moreover, given a complex-valued function $g$, we put $w_g := 1$ if $g$ is real, $w_g := \tfrac{1}{2}$ otherwise, and write

$$M(x; g) := \sum_{n \leqslant x} g(n), \qquad Z(x, g) := \sum_{p \leqslant x} \frac{g(p)}{p}.$$

**Theorem 1.3 ([10]).** *Let*

$$\mathfrak{a} \in \, ]0, \tfrac{1}{4}], \quad \mathfrak{b} \in [\mathfrak{a}, \tfrac{1}{2}[, \quad \mathfrak{h} := (1 - \mathfrak{b})/\mathfrak{b}, \quad A \geqslant 2\mathfrak{b}, \quad B > 0, \quad \beta := \beta_0(\mathfrak{b}, A),$$
$$x \geqslant 2, \quad 1/\sqrt{\log x} < \varepsilon \leqslant \tfrac{1}{2},$$

*and let the multiplicative functions $g$, $r$, such that $r \in \mathcal{M}(x; 2A, B)$, $|g| \leqslant r$, satisfy the conditions*

(1·12) $$\sum_{p \leqslant x} \frac{r(p) - \mathfrak{Re}\, g(p)}{p} \leqslant \tfrac{1}{2} \beta \mathfrak{b} \log(1/\varepsilon),$$

(1·13) $$\sum_{x^\varepsilon < p \leqslant y} \frac{\{r(p) - \mathfrak{Re}\, g(p)\}^{\mathfrak{h}} \log p}{p} \ll \varepsilon^{\delta \mathfrak{h}} \log y \qquad (x^\varepsilon < y \leqslant x),$$

*with $\delta \in [\mathfrak{a}, \tfrac{1}{3}\beta\mathfrak{b}]$, and*

(1·14) $$\min_{x^\varepsilon < p \leqslant x} r(p) \geqslant 4\mathfrak{b}.$$

*We then have*

(1·15) $$M(x; g) = M(x; r) \prod_p \frac{\sum_{p^\nu \leqslant x} g(p^\nu)/p^\nu}{\sum_{p^\nu \leqslant x} r(p^\nu)/p^\nu} + O\left(\frac{x\, \varepsilon^{w_g \delta}\, \mathrm{e}^{Z(x;r) - \mathfrak{c} Z(x;|g|-g)}}{\log x}\right)$$

*where $\mathfrak{c} := \mathfrak{b}/A$. The implicit constant in (1·15) depends at most upon $A$, $B$, $\mathfrak{a}$, and $\mathfrak{b}$.*

## 2. Proof of Theorem 1.1

As noted by Halász [5], we may assume that $f$ is integer-valued. (Note, however, that a slight modification of his construction is needed to ensure that changing the range of $f$ does not create new coincidences.) With this reduction, we plainly have

$$Q(x; f, \mu) \leqslant \int_{-1/2}^{1/2} |M(x; \vartheta)|\, \mathrm{d}\vartheta$$

with

$$M(x; \vartheta) := \sum_{n \leqslant x} \mu(n) \mathrm{e}^{2\pi i \vartheta f(n)}.$$

From Corollary III.4.12 in [8], we get, uniformly for $\vartheta \in \mathbb{R}$, $T \geqslant 1$, $x \geqslant 1$,

$$(2 \cdot 1) \qquad M(x; \vartheta) \ll \frac{x\{1 + m(x; \vartheta, T)\}}{e^{m(x; \vartheta, T)}} + \frac{x}{T},$$

where we have put

$$m(x; \vartheta, T) := \min_{|\tau| \leqslant T} \sum_{p \leqslant x} \frac{1 + \cos(2\pi\vartheta f(p) - \tau \log p)}{p}.$$

We select $T := \log x$, so that the second term on the right of (2·1) is negligible compared to the upper bound in (1·4). Let $h_\vartheta$ defined by

$$h_\vartheta(t) := 1 + \min\{\cos(t), \cos(2\pi\vartheta - t)\} \qquad (t \in \mathbb{R}),$$

so that

$$s_\vartheta := \frac{1}{2\pi} \int_{-\pi}^{\pi} h_\vartheta(t)\, \mathrm{d}t = 1 - \frac{2}{\pi} |\sin(\pi\vartheta)| \qquad (\vartheta \in [-\tfrac{1}{2}, \tfrac{1}{2}]),$$

and, for suitable $\tau \in [-T, T]$,

$$m(x; \vartheta, T) \geqslant \sum_{p \leqslant x} \frac{h_\vartheta(\tau \log p)}{p}.$$

The right-hand side may be estimated via partial summation as made explicit in lemma III.4.13 of [8]. For any $w \in [2, x]$ and $\vartheta \in [-\tfrac{1}{2}, \tfrac{1}{2}]$, we have

$$(2 \cdot 2) \qquad \sum_{w < p \leqslant x} \frac{h_\vartheta(\tau \log p)}{p} = s_\vartheta \log\left(\frac{\log x}{\log w}\right) + O\left(\frac{1}{w \log x} + \frac{1 + |\tau|}{e^{\sqrt{\log w}}}\right).$$

If $1 \leqslant |\tau| \leqslant T$, we select $w := (\log_2 x)^2$ to obtain

$$m(x; \vartheta, T) \geqslant s_\vartheta \log_2 x + O(\log_3 x).$$

Next, set

$$\log v := (\log x) \exp\left\{-\frac{2\cos^2(\pi\vartheta)E(x) + 2F(x)}{2 + s_\vartheta}\right\}.$$

If $1/\log v < |\tau| \leqslant 1$, we put $w := v$ in (2·2) and get

$$\sum_{v < p \leqslant x} \frac{h_\vartheta(\tau \log p)}{p} \geqslant \frac{2s_\vartheta \cos^2(\pi\vartheta)}{2 + s_\vartheta} E(x) + \frac{2s_\vartheta}{2 + s_\vartheta} F(x) + O(1).$$

And finally, if $|\tau| \leqslant 1/\log v$, we have trivially

$$\sum_{p \leqslant v} \frac{1 + \cos(2\pi\vartheta f(p) - \tau \log p)}{p} = \sum_{p \leqslant v} \frac{1 + \cos(2\pi\vartheta f(p))}{p} + O(1)$$

$$= (1 + \cos(2\pi\vartheta)) \sum_{\substack{p \leqslant v \\ f(p)=1}} \frac{1}{p} + 2 \sum_{\substack{p \leqslant v \\ f(p)=0}} \frac{1}{p} + O(1)$$

$$\geqslant 2\cos^2(\pi\vartheta) E(x) + 2F(x) - 2\log\left(\frac{\log x}{\log v}\right) + O(1)$$

$$\geqslant \frac{2s_\vartheta \cos^2(\pi\vartheta)}{2 + s_\vartheta} E(x) + \frac{2s_\vartheta}{2 + s_\vartheta} F(x) + O(1).$$

Therefore, we get in all cases

$$(2 \cdot 3) \qquad m(x; \vartheta, T) \geqslant \frac{2s_\vartheta \cos^2(\pi\vartheta)}{2 + s_\vartheta} E(x) + \frac{2s_\vartheta}{2 + s_\vartheta} F(x) + O(1)$$

$$\geqslant c\cos^2(\pi\vartheta) E(x) + cF(x) + O(1).$$

Integrating over $\vartheta$ immediately yields the result stated. $\square$

## 3. Proof of Theorem 1.2

Let us introduce the multiplicative function $g(n) := \mu(n)z^{f(n)}$ with $z := -\varrho e^{2\pi i\vartheta}$, $|\vartheta| \leqslant \frac{1}{2}$, $\kappa \leqslant \varrho \leqslant 1/\kappa$. Put $r(n) := \mu(n)^2 \varrho^{f(n)}$. From (2·3), we see that, with $c$ as in the statement of Theorem 1.1,

$$\sum_{p \leqslant x} \frac{r(p) - \Re e\,(g(p)/p^{i\tau})}{p} \geqslant c\varrho \sin^2(\pi\vartheta)E(x) + c\varrho F(x) + O(1) \qquad (|\tau| \leqslant T := \log x).$$

We may therefore apply Corollary 2.1 of [10] to get

$$(3\cdot 1) \qquad M(x;g) \ll M(x;r)\left\{e^{-c\varrho E(x)\sin^2(\pi\vartheta) - c\varrho F(x)}\log_2 x + \frac{1}{(\log x)^\kappa}\right\}.$$

With the aim of applying Cauchy's formula to detect $N_m(x; f, \mu)$, we next seek an estimate for $M(x;g)$ when $\vartheta$ is small, namely

$$|\vartheta| \leqslant \vartheta_0 := K\sqrt{\frac{\log_3 x}{\log_2 x}},$$

where $K$ is a large constant—actually any $K > 1/\sqrt{4\kappa c}$ will do. We have

$$\sum_{p \leqslant x} \frac{r(p) - \Re e\,g(p)}{p} = \varrho(1 - \cos 2\pi\vartheta)E(x) + 2\varrho F(x) \leqslant 2\varrho\pi^2\vartheta^2 + 2\varrho F(x),$$

hence condition (1·12) is plainly fulfilled with $\varepsilon := |\vartheta|^{2/\delta} + (\log_2 x)^{-c_0/(\mathfrak{h}\delta)}$ provided $\delta$ is chosen sufficiently small in terms of $\mathfrak{b}$, $\kappa$ and $K$. Next, for $x^\varepsilon < y \leqslant x$, we have

$$\sum_{x^\varepsilon < p \leqslant y} \frac{\{r(p) - \Re e\,g(p)\}^\mathfrak{h} \log p}{p} \ll \varrho\vartheta^{2\mathfrak{h}}\log y + \varrho \sum_{\substack{x^\varepsilon < p \leqslant y \\ f(p) = 0}} \frac{\log p}{p}$$

$$\ll_\kappa \left\{\varepsilon^{\delta\mathfrak{h}} + (\log_2 x)^{-c_0}\right\}\log y,$$

so hypothesis (1·13) is also verified. Since (1·14) holds trivially on selecting $\mathfrak{b} := \kappa/4$, and hence $\mathfrak{h} = 4/\kappa - 1$, we conclude that (1·15) is valid. We obtain, with $\mathfrak{c} := \kappa\mathfrak{b}$,

$$(3\cdot 2) \quad M(x;g) = M(x;r) \prod_{\substack{p \leqslant x \\ f(p) = 1}} \frac{1 - z/p}{1 + \varrho/p} \prod_{\substack{p \leqslant x \\ f(p) = 0}} \frac{1 - 1/p}{1 + 1/p} + O\left(\frac{x\varepsilon^{\delta/2}e^{Z(x;r) - \mathfrak{c}Z(x;r-g)}}{\log x}\right).$$

Now, appealing for instance to theorem 1.1 of [9], we observe that

$$M(x;r) \asymp \frac{xe^{Z(x;r)}}{\log x}$$

and so we may rewrite (3·2) as

$$M(x;g) = M(x;r)\left\{\lambda_f e^{-(z+\varrho)E(x) - 2F(x)} + O\left(\left(|\vartheta| + (\log_2 x)^{-c_0\mathfrak{h}/2}\right)e^{-c_1\vartheta^2 E(x) - c_1 F(x)}\right)\right\},$$

valid for $|\vartheta| \leqslant \vartheta_0$ and some constant $c_1 > 0$. Integrating on the circle $|z| = \varrho := m/E(x)$ and taking (3·1) into account, we readily obtain in the stated range for $m$,

$$\begin{aligned}(3\cdot 3) \quad N_m(x;f,\mu) &= (-1)^m \int_{-1/2}^{1/2} e^{-2i\pi\vartheta m} \varrho^{-m} M(x;g)\, d\vartheta \\ &= (-1)^m \lambda_f M(x;r) \frac{E(x)^m}{m!\, e^m} \left\{ e^{-2F(x)} + O\!\left(\frac{e^{-c_2 F(x)}}{(\log_2 x)^b}\right)\right\},\end{aligned}$$

with $c_2 := \min(c_1, c\kappa)$. Since, by a straightforward variant of corollary 2.4 of [10][2],

$$N_m(x;f) = M(x;r)\frac{E(x)^m}{m!\, e^m}\left\{1 + O\!\left(\frac{1}{\sqrt{\log_2 x}}\right)\right\},$$

we reach the required conclusion.

Régis de la Bretèche  
Institut de Mathématiques de Jussieu  
UMR 7586  
Université Paris Diderot-Paris 7  
Sorbonne Paris Cité,  
Case 7012, F-75013 Paris  
France

Gérald Tenenbaum  
Institut Élie Cartan  
Université de Lorraine  
BP 70239  
54506 Vandœuvre-lès-Nancy Cedex  
France

---

2. Applied to $\omega(n;E)$ instead of $\Omega(n;E)$ with the notation of [10].